\documentclass[11pt]{article}

\usepackage[T1]{fontenc}
\usepackage{lmodern}
\usepackage{amsmath,amssymb,amsthm,mathtools,bm}
\usepackage[a4paper,margin=25mm]{geometry}
\usepackage{enumitem}
\usepackage{microtype}
\usepackage{hyperref}
\hypersetup{colorlinks=true,linkcolor=blue,citecolor=blue,urlcolor=blue}

\newtheorem{theorem}{Theorem}[section]
\newtheorem{lemma}[theorem]{Lemma}
\newtheorem{proposition}[theorem]{Proposition}
\newtheorem{corollary}[theorem]{Corollary}
\theoremstyle{definition}

\newtheorem{remark}[theorem]{Remark}

\newcommand{\R}{\mathbb{R}}
\newcommand{\ip}[2]{\left\langle #1,#2\right\rangle}
\newcommand{\1}{\mathbf{1}}
\newcommand{\diag}{\operatorname{diag}}

\newcommand{\dist}{\operatorname{dist}}

\title{Least-distortion Euclidean embeddings of distance-regular graphs}
\author{Hirotake Kurihara\thanks{Department of Applied Science, Yamaguchi University,
2-16-1 Tokiwadai, Ube 755-8611, Japan.
Email: \href{mailto:kurihara-hiro@yamaguchi-u.ac.jp}{\texttt{kurihara-hiro@yamaguchi-u.ac.jp}}.}}
\date{}

\begin{document}
\maketitle

\begin{abstract}
Let $\Gamma$ be a distance-regular graph of diameter $d$, with
distinct eigenvalues
$k=\theta_0>\theta_1>\cdots>\theta_d$ and cosine sequences
$u_i(\theta_j)$. For every $1\le i\le d$, we prove that
$(1-u_i(\theta_j))/(1-u_1(\theta_j))$ is maximized over
$1\le j\le d$ at $j=1$.
This proves a conjecture of Cioab\u{a}, Gupta, Ihringer, and Kurihara
arising from Vallentin's work on Euclidean embeddings.
As a consequence, we determine the least Euclidean distortion of
$\Gamma$ and prove that the standard representation associated with
$\theta_1$ is optimal.
The proof uses the Perron--Frobenius theorem and a nonnegative matrix
representation of an element of the Bose--Mesner algebra.
\end{abstract}



\section{Introduction}
\label{sec:introduction}

In this paper, we determine a least-distortion Euclidean embedding of
every distance-regular graph. We prove that the standard representation
associated with the second largest adjacency eigenvalue is optimal.
We also give an exact formula for the least distortion in terms of the
corresponding cosine sequence.

The least Euclidean distortion of a finite metric space measures how
closely its distances can be represented by Euclidean distances.
Enflo~\cite{Enflo} proved that the least Euclidean distortion of the
$n$-dimensional hypercube is $\sqrt{n}$.
Bourgain~\cite{Bourgain} proved that every metric space on $n$ points
admits a Euclidean embedding with distortion $O(\log n)$.
For graph metrics, Linial, London, and Rabinovich~\cite{LLR}
developed connections with combinatorial optimization and expressed
the least Euclidean distortion as a semidefinite program.
Linial and Magen~\cite{LM} studied optimal embeddings of products
of cycles and expanders. Their results for expanders show that
Bourgain's bound has the correct order of growth.

Lower bounds for distortion also reflect combinatorial and spectral
properties of graphs. Linial, Magen, and Naor~\cite{LMN} proved
a lower bound of order $\sqrt{g}$ for the Euclidean distortion of
$k$-regular graphs of girth $g$, where $k\ge3$.
Grigorchuk and Nowak~\cite{GN} related the $\ell_p$-distortion
to the diameter and the first positive eigenvalue of the discrete
$p$-Laplacian. Jolissaint and Valette~\cite{JV} obtained a lower
bound in terms of the $p$-spectral gap and the maximal displacement
of vertex permutations. For $p=2$, their bound gives the exact
distortion of even cycles and hypercubes.

Distance-regular graphs provide a setting in which this optimization
problem can be studied through adjacency eigenvalues and cosine
sequences. Vallentin~\cite{Vallentin} used the Bose--Mesner algebra
to obtain a spectral lower bound for their least Euclidean distortion.
He also determined the distortion of Hamming graphs, Johnson graphs,
and connected strongly regular graphs.
Kobayashi and Kondo~\cite{KK} obtained exact values for generalized
polygons by semidefinite programming.
Vallentin's proposed equality for the general lower bound does not hold
in general. Cioab\u{a}, Gupta, Ihringer, and Kurihara~\cite{CGIK}
constructed examples of diameter $d$ for which an optimal embedding
has its largest contraction at distance $d-1$.
They also proposed revised conjectures on the cosine sequences.

We prove their conjecture that identifies the eigenvalue maximizing a
ratio of cosine differences.
Let $\Gamma$ be a finite connected distance-regular graph of diameter
$d\ge1$ and valency $k$.
Write its distinct adjacency eigenvalues as
$k=\theta_0>\theta_1>\cdots>\theta_d$.
Let $A_i=v_i(A_1)$ be its distance-$i$ matrix, let $k_i=v_i(k)$,
and put $u_i(\theta_j)=v_i(\theta_j)/k_i$.
The sequence
$u_0(\theta_j),u_1(\theta_j),\ldots,u_d(\theta_j)$ is called the
\emph{cosine sequence} for $\theta_j$.
These definitions are recalled in Section~\ref{sec:preliminaries}.
In particular, $u_0(\theta_j)=1$ and $u_1(\theta_j)=\theta_j/k$.

For $1\le i,j\le d$, set
\begin{equation}
\label{eq:S}
  S_i(\theta_j):=\frac{1-u_i(\theta_j)}{1-u_1(\theta_j)}.
\end{equation}
The denominator is positive because $\theta_j<k$.

The following Theorem~\ref{thm:central} is the main result of this paper.
This theorem resolves Conjecture~1.7 in~\cite{CGIK} affirmatively.

\begin{theorem}
\label{thm:central}
Let $\Gamma$ be a finite connected distance-regular graph of
diameter $d\ge1$. Then, for all $1\le i,j\le d$,
\[
  S_i(\theta_j)\le S_i(\theta_1).
\]
Moreover, equality holds if and only if $i=1$ or $j=1$.
\end{theorem}

We prove Theorem~\ref{thm:central} in Section~\ref{sec:proof}.
The proof uses the Bose--Mesner algebra.
For a fixed vertex, we consider the space spanned by the
characteristic vectors of its distance layers.
We construct a basis for the orthogonal complement of the
constant vector in this space.
In this basis, certain elements of the Bose--Mesner algebra have
nonnegative matrix representations with eigenvalues $S_i(\theta_j)$.
The differences of the cosine sequence for $\theta_1$ give a
positive eigenvector. The Perron--Frobenius theorem identifies its
eigenvalue as the spectral radius and gives the equality condition.
As a corollary, we obtain
$k+\theta_1+\theta_d\ge a_1+a_2$ for $d\ge2$.
Equality holds if and only if $d=2$.

For an injective map $F$ from the vertex set of $\Gamma$ to a
Euclidean space, define
\[
  \dist F=
  \left(\max_{x\ne y}\frac{\|F(x)-F(y)\|}{\partial(x,y)}\right)
  \left(\max_{x\ne y}\frac{\partial(x,y)}{\|F(x)-F(y)\|}\right),
\]
where $\partial$ is the graph distance. The \emph{least Euclidean
distortion} $c_2(\Gamma)$ is the infimum of $\dist F$ over all such
maps and all finite-dimensional Euclidean target spaces.

The standard representation $\iota_j$ associated with $\theta_j$ maps each
vertex to the orthogonal projection of its coordinate vector onto the
$\theta_j$-eigenspace of $A_1$, normalized to have unit length.
The inner product of the images of two vertices at distance $i$
is $u_i(\theta_j)$.
The standard representation associated with $\theta_1$ is injective.

As a corollary of Theorem~\ref{thm:central}, we obtain the following
formula for the least Euclidean distortion.

\begin{corollary}
\label{cor:distortion}
Let $\Gamma$ be as in Theorem~\ref{thm:central}. Then
\begin{equation}
\label{eq:distortion}
  c_2(\Gamma)^2
  =\max_{\lfloor d/2\rfloor+1\le i\le d}
    i^2\frac{1-u_1(\theta_1)}{1-u_i(\theta_1)}.
\end{equation}
The standard representation associated with $\theta_1$ attains this
distortion.
\end{corollary}

In~\eqref{eq:distortion}, it suffices to take the maximum over
distances $i>d/2$.
Conjecture~1.6 in~\cite{CGIK} asserts that the maximum is attained
at $i=d-1$ or $i=d$.
Our result establishes the optimality of the standard representation
associated with $\theta_1$, but does not settle
Conjecture~1.6 in~\cite{CGIK}.

Optimality does not imply uniqueness up to similarity.
At the end of Section~\ref{sec:distortion}, we construct an optimal
embedding of the Sylvester graph in $\R^8$.
The affine span of its standard representation associated with
$\theta_1$ has dimension $16$.

The classical Wirtinger inequality bounds the $L^2$-norm of a
mean-zero periodic function by the $L^2$-norm of its derivative.
Fan, Taussky, and Todd~\cite{FTT} established sharp discrete analogues
for finite sequences. These analogues can also be expressed as
inequalities between the squared side lengths and squared diagonal
lengths of a closed Euclidean polygon; see Pech~\cite{Pech}.
We use a version for the distance relations of $\Gamma$.
For a nonconstant map $F:X\to\R^m$, let $W_i(F)$ denote the mean
squared distance between the images of ordered pairs at graph distance
$i$.
In Section~\ref{sec:distortion}, for every $1\le i\le d$, we prove
the Wirtinger inequality
\[
  \frac{W_i(F)}{W_1(F)}
  \le S_i(\theta_1)
  =\frac{W_i(\iota_1)}{W_1(\iota_1)}.
\]
Thus the standard representation $\iota_1$ realizes the sharp constant.
For an embedding of distortion $D$, this Wirtinger inequality gives
$D^2\ge i^2/S_i(\theta_1)$.
The case $i=d$ recovers Vallentin's lower bound~\cite{Vallentin},
and allowing all $i$ gives its extension in~\cite{CGIK}.
Together with the direct calculation of $\dist\iota_1$, this proves
the distortion formula in Corollary~\ref{cor:distortion}.

\section{Preliminaries}
\label{sec:preliminaries}

All graphs in this paper are finite, simple,
undirected, and connected. We assume that they have at least two
vertices.

We follow the standard terminology for distance-regular graphs
in~\cite{BCN,Godsil,VDKT}.
A graph $\Gamma=(X,E)$ of diameter $d$ is \emph{distance-regular}
if, whenever $\partial(x,y)=i$, the numbers of neighbors of $y$
at distances $i-1$, $i$, and $i+1$ from $x$ depend only on $i$.
These numbers are denoted by $c_i$, $a_i$, and $b_i$, respectively.
We use the conventions $c_0=0$ and $b_d=0$.
Then $b_0=k$ is the valency, $c_1=1$, and
$a_i+b_i+c_i=k$ for $0\le i\le d$.
Moreover, $b_i>0$ for $0\le i<d$, and $c_i>0$ for $1\le i\le d$.
The array
$\{b_0,\ldots,b_{d-1};c_1,\ldots,c_d\}$ is called the
\emph{intersection array} of $\Gamma$.

For $x\in X$, let
\[
  \Gamma_i(x)=\{y\in X:\partial(x,y)=i\},\qquad k_i=|\Gamma_i(x)|.
\]
The numbers $k_i$ are independent of $x$.
It is known that we have
\begin{equation}
\label{eq:ki}
  k_0=1,\qquad k_1=k,\qquad
  k_{i-1}b_{i-1}=k_i c_i\quad(1\le i\le d).
\end{equation}
More generally, the intersection numbers of $\Gamma$ are given by
\[
  p_{ij}^{h}
  =|\{z\in X:\partial(x,z)=i,\ \partial(y,z)=j\}|
  \quad\text{when }\partial(x,y)=h.
\]
Distance-regularity implies that they are independent of the choice
of $x$ and $y$. In this notation,
\[
  b_i=p_{i+1,1}^{i},\qquad a_i=p_{i,1}^{i},\qquad c_i=p_{i-1,1}^{i},
\]
with the boundary terms interpreted as zero.
We also use the standard identity
\begin{equation}
\label{lem:reversible}
  k_i p_{jh}^{i}=k_j p_{ih}^{j}\qquad(0\le i,j,h\le d).
\end{equation}
For $0\le i\le d$, the \emph{distance-$i$ matrix} is defined by
\[
  (A_i)_{xy}=\begin{cases}
    1,&\partial(x,y)=i,\\
    0,&\partial(x,y)\ne i.
  \end{cases}
\]
Thus $A_0=I$, $A_1$ is the adjacency matrix, and
$\sum_{i=0}^d A_i=J$, where $I$ is the identity matrix and
$J$ is the all-ones matrix.
The distance matrices satisfy
\[
  A_iA_j=\sum_{h=0}^d p_{ij}^{h}A_h.
\]
Then the span $\mathfrak{A}=\operatorname{span}_{\R}\{A_0,\ldots,A_d\}
\subseteq\operatorname{Mat}_X(\R)$
is a commutative algebra, called the \emph{Bose--Mesner algebra}.
The recurrence
\begin{equation}
\label{3-term}
  A_1A_i=b_{i-1}A_{i-1}+a_iA_i+c_{i+1}A_{i+1}
\end{equation}
holds for $0\le i\le d$, with terms outside this range omitted.
It gives polynomials $v_i$ of degree $i$ such that $A_i=v_i(A_1)$.
In particular, $v_0(t)=1$, $v_1(t)=t$, and $v_i(k)=k_i$.
Thus $\mathfrak{A}=\R[A_1]$.

The matrix $A_1$ has exactly $d+1$ distinct eigenvalues
$k=\theta_0>\theta_1>\cdots>\theta_d$.
Let $E_j$ be the orthogonal projection onto its $\theta_j$-eigenspace,
and put $m_j=\operatorname{rank}E_j$.
The matrices $E_j$ are the primitive idempotents of the
Bose--Mesner algebra. They satisfy
\[
  E_iE_j=\delta_{ij}E_j,\qquad
  \sum_{j=0}^d E_j=I,\qquad E_0=\frac1{|X|} J.
\]
Here $\delta_{ij}$ is the Kronecker delta.
The \emph{cosine sequence}, also called the \emph{standard sequence},
for $\theta_j$ is defined by
\[
  u_i(\theta_j)=\frac{v_i(\theta_j)}{k_i}\qquad(0\le i\le d).
\]
The spectral expansions are
\begin{equation}
\label{eq:idempotents}
  A_i=\sum_{j=0}^d k_i u_i(\theta_j)E_j,
  \qquad
  E_j=\frac{m_j}{|X|}\sum_{h=0}^d u_h(\theta_j)A_h.
\end{equation}
In particular, $A_iE_j=k_i u_i(\theta_j)E_j$.
The cosine sequences satisfy
\begin{equation}
\label{eq:cosine-recurrence}
  \theta_j u_i(\theta_j)
  =c_i u_{i-1}(\theta_j)+a_i u_i(\theta_j)+b_i u_{i+1}(\theta_j),
\end{equation}
where the first term is omitted for $i=0$ and the last for $i=d$.

For $1\le i,j\le d$, put
$h_i^{(j)}=u_{i-1}(\theta_j)-u_i(\theta_j)$ and
$h^{(j)}=(h_1^{(j)},\ldots,h_d^{(j)})^T$.
We use the following standard fact; see Godsil~\cite[Lemma~13.2.1]{Godsil}.
\begin{lemma}
\label{lem:monotone}
Every entry of $h^{(1)}$ is positive. Equivalently,
\[
  1=u_0(\theta_1)>u_1(\theta_1)>\cdots>u_d(\theta_1).
\]
\end{lemma}

Let $\R^X$ be the standard vector space indexed by $X$ over $\R$,
with its usual inner product $\ip{\cdot}{\cdot}$.
Let $\1$ be the all-ones vector in $\R^X$, and let $e_x$ be the
coordinate vector indexed by $x\in X$.
An $\mathfrak{A}$-module is a subspace of $\R^X$ invariant under the left action of $\mathfrak{A}$.
Fix a vertex $x\in X$.

Write
\[
  s_i=A_i e_x\qquad(0\le i\le d)
\]
for the characteristic vector of $\Gamma_i(x)$.
Then $\ip{s_i}{s_j}=k_i\delta_{ij}$ and $\sum_{i=0}^d s_i=\1$.
We use the $\mathfrak{A}$-module
\begin{equation}
\label{eq:primary}
  W=\mathfrak{A} e_x
   =\operatorname{span}_{\R}\{s_0,\ldots,s_d\}.
\end{equation}
We also use the $\mathfrak{A}$-module
\begin{equation}
\label{eq:W0}
  W_0=(I-E_0)W=W\cap\1^\perp
      =\bigoplus_{j=1}^d E_jW.
\end{equation}
Each $E_jW$ is spanned by $E_je_x$.
Since $E_1e_x,\ldots,E_de_x$ are linearly independent,
we have $\dim W_0=d$.
In Section~\ref{sec:proof}, we construct a basis $z_1,\ldots,z_d$
of $W_0$ and restrict elements of $\mathfrak{A}$ to $W_0$.

\section{Proof of the main theorem}
\label{sec:proof}

Since $A_1$ is symmetric and $A_1\1=k\1$, every $w\in W$ satisfies
\[
  \ip{\1}{\left(I-\frac{A_1}{k}\right)w}
  =\ip{\left(I-\frac{A_1}{k}\right)\1}{w}=0.
\]
It follows that $(I-A_1/k)W\subseteq W_0$.
For $1\le p\le d$, put
\begin{equation}
\label{eq:eta}
  \eta_p=\sum_{i=0}^{p-1}s_i,\qquad
  z_p=\left(I-\frac{A_1}{k}\right)\eta_p.
\end{equation}
Thus $\eta_p$ is the characteristic vector of the vertices at
distance less than $p$ from $x$, and $z_p\in W_0$.
Equation~\eqref{3-term} gives
$A_1s_i=b_{i-1}s_{i-1}+a_is_i+c_{i+1}s_{i+1}$.
Using $a_i+b_i+c_i=k$, we obtain
\begin{align*}
  z_p
  &=\frac1k\sum_{i=0}^{p-1}(ks_i-A_1s_i)\\
  &=\frac1k\left(
      \sum_{i=0}^{p-2}(k-b_i-a_i-c_i)s_i
      +(k-a_{p-1}-c_{p-1})s_{p-1}-c_ps_p\right)\\
  &=\frac1k\left(b_{p-1}s_{p-1}-c_ps_p\right).
\end{align*}
By~\eqref{eq:ki}, the numbers
$C_p=\frac{k_{p-1}b_{p-1}}{k}=\frac{k_pc_p}{k}$
are positive for $1\le p\le d$. Hence we have
\begin{equation}
\label{eq:z}
  z_p=C_p\left(\frac{s_{p-1}}{k_{p-1}}-\frac{s_p}{k_p}\right).
\end{equation}
Put $C=\diag(C_1,\ldots,C_d)$.

\begin{lemma}
\label{lem:basis}
For $1\le p,q\le d$,
\begin{equation}
\label{eq:dual-basis}
  \ip{\eta_p}{z_q}=C_q\delta_{pq}.
\end{equation}
Moreover, $z_1,\ldots,z_d$ form a basis of $W_0$.
If $w=\sum_{q=1}^d\alpha_qz_q\in W_0$, then
\[
  \alpha_p=\frac{\ip{\eta_p}{w}}{C_p}.
\]
\end{lemma}

\begin{proof}
Since $\ip{s_i}{s_j}=k_i\delta_{ij}$,
\eqref{eq:eta} and~\eqref{eq:z} give
\[
  \ip{\eta_p}{z_q}
  =C_q\sum_{i=0}^{p-1}
       \left(\delta_{i,q-1}-\delta_{i,q}\right)
  =C_q\delta_{pq}.
\]
As $C_q>0$, this identity implies that $z_1,\ldots,z_d$
are linearly independent. They form a basis because $\dim W_0=d$.
The coordinate formula also follows from~\eqref{eq:dual-basis}.
\end{proof}

For $1\le j\le d$, put
\[
  \psi_j=\frac{|X|}{m_j}E_j e_x
        =\left(\sum_{i=0}^d u_i(\theta_j)A_i\right)e_x
        =\sum_{i=0}^d u_i(\theta_j)s_i.
\]
Since $E_0E_j=0$ for $1\le j\le d$, the vectors
$\psi_1,\ldots,\psi_d$ belong to $W_0$.
They form a basis of $W_0$.
By the definitions of $h_p^{(j)}$ and $\eta_p$,
\[
  \psi_j=u_d(\theta_j)\1+\sum_{p=1}^d h_p^{(j)}\eta_p.
\]
Applying $I-A_1/k$ gives
\begin{equation}
\label{eq:cosine-coordinates}
  \left(1-\frac{\theta_j}{k}\right)\psi_j
  =\sum_{p=1}^d h_p^{(j)}z_p.
\end{equation}
Since $1-\theta_j/k>0$, the vectors
$h^{(1)},\ldots,h^{(d)}$ form a basis of $\R^d$.

For $1\le i\le d$, define
\begin{equation}
\label{eq:sigma}
  \sigma_i(t)=\frac{1-v_i(t)/k_i}{1-t/k}.
\end{equation}
Since $v_i(k)=k_i$, cancellation of the factor $t-k$ extends
$\sigma_i$ to a polynomial of degree $i-1$.
Consider the element
\[
  F_i=\sigma_i(A_1)\in \mathfrak{A}.
\]
Equation~\eqref{eq:sigma} gives
\begin{equation}
\label{eq:factor}
  F_i\left(I-\frac{A_1}{k}\right)=I-\frac{A_i}{k_i}.
\end{equation}
Also,
\begin{equation}
\label{eq:F-spectrum}
  F_i\psi_j=S_i(\theta_j)\psi_j
  \qquad(1\le j\le d).
\end{equation}
Thus Theorem~\ref{thm:central} compares the eigenvalues of $F_i|_{W_0}$.
We now give a nonnegative matrix representation of this restriction.

For $1\le i\le d$, define the symmetric matrix $M_i$ of size $d$ by
\begin{equation}
\label{eq:M}
  (M_i)_{pq}=
  \sum_{\substack{0\le h<\ell\le d\\
                   h<\min\{p,q\}\\\ell\ge\max\{p,q\}}}
    \frac{k_h p_{\ell i}^{h}}{k_i}
  \qquad(1\le p,q\le d).
\end{equation}
\begin{proposition}
\label{prop:nonnegative}
Let $1\le i\le d$.
The representation matrix of $F_i|_{W_0}$ with respect to $z_1,\ldots,z_d$
is $C^{-1}M_i$.
With respect to the basis
\[
  \frac{z_1}{\sqrt{C_1}},\ldots,\frac{z_d}{\sqrt{C_d}},
\]
its representation matrix is
\begin{equation}
\label{eq:N}
  N_i=C^{-1/2}M_iC^{-1/2}.
\end{equation}
In particular, $N_i$ is symmetric and entrywise nonnegative.
\end{proposition}

\begin{proof}
By~\eqref{eq:eta} and~\eqref{eq:factor}, we have
\[
  F_i z_q=\left(I-\frac{A_i}{k_i}\right)\eta_q.
\]
Lemma~\ref{lem:basis} shows that the $p$th coordinate of $F_i z_q$
in the basis $z_1,\ldots,z_d$ is
\begin{equation}
\label{eq:matrix-entry}
  \frac1{C_p}
  \ip{\eta_p}{\left(I-\frac{A_i}{k_i}\right)\eta_q}.
\end{equation}
Suppose first that $p\le q$.
On the support of $\eta_p$, both $\eta_p$ and $\eta_q$ have value one.
Since every row of $A_i$ sums to $k_i$, we obtain
\begin{align}
\label{eq:crossing}
  \ip{\eta_p}{\left(I-\frac{A_i}{k_i}\right)\eta_q}
  &=\sum_{\substack{y\in X\\\partial(x,y)<p}}
    \left(1-\frac1{k_i}
      \sum_{\substack{z\in X\\\partial(x,z)<q}}(A_i)_{yz}\right)
      \notag\\
  &=\frac1{k_i}\sum_{\substack{y\in X\\\partial(x,y)<p}}
    \left(\sum_{z\in X}(A_i)_{yz}
      -\sum_{\substack{z\in X\\\partial(x,z)<q}}(A_i)_{yz}\right)
      \notag\\
  &=\frac1{k_i}
    \sum_{\substack{y,z\in X\\
                    \partial(x,y)<p,\ \partial(x,z)\ge q}}
      (A_i)_{yz}\notag\\
  &=\sum_{h=0}^{p-1}\sum_{\ell=q}^d\frac{k_h p_{\ell i}^{h}}{k_i}
   =(M_i)_{pq}.
\end{align}
For the fourth equality, fix $y\in\Gamma_h(x)$.
There are $p_{\ell i}^{h}$ vertices $z\in\Gamma_\ell(x)$ at distance $i$
from $y$, and there are $k_h$ choices for $y$.

For $p>q$, the same identity follows by symmetry:
\[
  \ip{\eta_p}{\left(I-\frac{A_i}{k_i}\right)\eta_q}
  =\ip{\eta_q}{\left(I-\frac{A_i}{k_i}\right)\eta_p}
  =(M_i)_{qp}=(M_i)_{pq}.
\]
Thus~\eqref{eq:matrix-entry} is $(M_i)_{pq}/C_p$.
Rescaling each basis vector $z_p$ by $1/\sqrt{C_p}$
gives~\eqref{eq:N}.
Every summand in~\eqref{eq:M} is nonnegative, so $N_i$
is entrywise nonnegative.
\end{proof}

\begin{proof}[Proof of Theorem~\ref{thm:central}]
For $i=1$, we have $S_1(\theta_j)=1$ for every $j$.
Suppose that $2\le i\le d$.
By Proposition~\ref{prop:nonnegative}, the matrix $N_i$
is symmetric and entrywise nonnegative.

We first show that $N_i$ is irreducible.
Choose a geodesic $x_0,x_1,\ldots,x_d$ in $\Gamma$.
Fix $1\le p<d$ and put $h=\max\{0,p+1-i\}$.
Since $i\ge2$, we have $0\le h\le p-1$ and $p+1\le h+i\le d$.
The distances from $x_0$ to $x_h$ and $x_{h+i}$ are $h$ and $h+i$,
respectively, and $\partial(x_h,x_{h+i})=i$.
Thus $p_{h+i,i}^{h}>0$, and~\eqref{eq:M} gives
\[
  (M_i)_{p,p+1}\ge\frac{k_h p_{h+i,i}^{h}}{k_i}>0.
\]
The positive superdiagonal entries and the symmetry of $M_i$ show
that $M_i$ is irreducible.
Positive diagonal scaling preserves the positions of its nonzero entries.
Hence $N_i$ is also irreducible.

Equation~\eqref{eq:cosine-coordinates} gives
\[
  \psi_j=\frac1{1-\theta_j/k}
    \sum_{p=1}^d\sqrt{C_p}\,h_p^{(j)}\frac{z_p}{\sqrt{C_p}}.
\]
Thus the coordinate vector of $\psi_j$ in the basis
$z_1/\sqrt{C_1},\ldots,z_d/\sqrt{C_d}$ is
$(1-\theta_j/k)^{-1}C^{1/2}h^{(j)}$.
By~\eqref{eq:F-spectrum} and Proposition~\ref{prop:nonnegative}, we have
\[
  N_i(C^{1/2}h^{(j)})
  =S_i(\theta_j)(C^{1/2}h^{(j)})
  \qquad(1\le j\le d).
\]
Since $N_i$ has size $d$ and
$C^{1/2}h^{(1)},\ldots,C^{1/2}h^{(d)}$ form a basis of $\R^d$,
the eigenvalues of $N_i$ are $S_i(\theta_1),\ldots,S_i(\theta_d)$.
Lemma~\ref{lem:monotone} shows that $C^{1/2}h^{(1)}$
is a positive eigenvector with eigenvalue $S_i(\theta_1)$.

By the Perron--Frobenius theorem, $S_i(\theta_1)$ is the spectral
radius of $N_i$ and is a simple eigenvalue.
Therefore,
\[
  S_i(\theta_j)<S_i(\theta_1)\qquad(2\le j\le d).
\]
This proves the inequality and its equality condition.
\end{proof}

\begin{remark}
\label{rem:probabilistic}
The proof also has a probabilistic interpretation.
The matrix $P_i=A_i/k_i$ is the Markov operator for the random walk
that moves uniformly to a vertex at distance $i$.
The $\mathfrak{A}$-module $W$ consists of the functions constant on
each $\Gamma_h(x)$, so it is the radial space about $x$.
The restriction to $W$ of the form
\[
  \mathcal Q_i(f,g)=\ip{f}{(I-P_i)g}
\]
is the radial Dirichlet form for the counting measure.
When functions are specified by their values on the distance
layers, this measure gives the weights $k_h$.
The calculation in Proposition~\ref{prop:nonnegative} gives
$(M_i)_{pq}=\mathcal Q_i(\eta_p,\eta_q)$.
Theorem~\ref{thm:central} gives the comparison
\[
  \mathcal Q_i(f,f)\le S_i(\theta_1)\mathcal Q_1(f,f)
  \qquad(f\in W,\ 1\le i\le d).
\]
\end{remark}

Theorem~\ref{thm:central} also gives the following inequality
for distance-regular graphs.
Since the intersection array determines the eigenvalues, this inequality
gives a necessary condition for its realization by a distance-regular graph.

\begin{corollary}
\label{cor:kurihara}
Let $\Gamma$ be a distance-regular graph of diameter $d\ge2$.
Then
\begin{equation}
\label{eq:kurihara}
  k+\theta_1+\theta_d\ge a_1+a_2.
\end{equation}
Equality holds if and only if $d=2$.
\end{corollary}

\begin{proof}
Suppose first that $d\ge3$.
The recurrence for the distance polynomials gives
\[
  v_2(t)=\frac{t^2-a_1t-k}{c_2},\qquad
  v_3(t)=\frac{(t-a_2)(t^2-a_1t-k)-b_1c_2t}{c_2c_3}.
\]
Since $k_3=kb_1b_2/(c_2c_3)$, each nontrivial eigenvalue $t$ satisfies
\[
  S_3(t)=
  \frac{t^2+(k-a_1-a_2)t+k^2-k(a_1+a_2)+a_1a_2-k-b_1c_2}
       {b_1b_2}.
\]
Consequently,
\begin{equation}
\label{eq:S3-difference}
  S_3(\theta_1)-S_3(\theta_d)
  =\frac{(\theta_1-\theta_d)
    (k+\theta_1+\theta_d-a_1-a_2)}{b_1b_2}.
\end{equation}
The left-hand side is positive by Theorem~\ref{thm:central},
since $d\ge3$.
As $\theta_1>\theta_d$ and $b_1b_2>0$, the inequality
in~\eqref{eq:kurihara} is strict.
For $d=2$, equality follows from the standard identity
\[
  \sum_{j=0}^d\theta_j=\sum_{i=1}^d a_i.
\]
\end{proof}

\section{Least Euclidean distortion}
\label{sec:distortion}

We first derive Corollary~\ref{cor:distortion} from
Theorem~\ref{thm:central}.
The \emph{standard representation} $\iota_j$ associated with
$\theta_j$ is
\[
  \iota_j:x\longmapsto \sqrt{\frac{|X|}{m_j}}\,E_je_x.
\]
Its image lies on the unit sphere in $E_j\R^X$, and the inner
product of the images of $x$ and $y$ is
$u_i(\theta_j)$ whenever $\partial(x,y)=i$.
A standard representation need not be injective.
However, $\iota_1$ is injective by Lemma~\ref{lem:monotone}.
If $\partial(x,y)=i$, then
\eqref{eq:idempotents} gives
\begin{equation}
\label{eq:standard-distances}
  \|\iota_1(x)-\iota_1(y)\|^2
  =2(1-u_i(\theta_1)).
\end{equation}
Thus adjacent vertices are mapped to points at distance
$\sqrt{2(1-u_1(\theta_1))}$.

The following lemma restates Lemma~2.1 in
\cite{CGIK} in our notation.
\begin{lemma}
\label{lem:distance-range}
Let $p$ and $q$ be positive integers such that $pq\le d$. Then
\[
  \frac{\sqrt{2(1-u_{pq}(\theta_1))}}{pq}
  \le \frac{\sqrt{2(1-u_q(\theta_1))}}{q}.
\]
Moreover, for $2\le i\le d$,
\[
  \frac{\sqrt{2(1-u_i(\theta_1))}}{i}
  < \sqrt{2(1-u_1(\theta_1))}.
\]
\end{lemma}
\begin{proof}
Choose a geodesic $x_0,x_1,\ldots,x_{pq}$ in $\Gamma$.
The triangle inequality and~\eqref{eq:standard-distances} give
\[
  \sqrt{2(1-u_{pq}(\theta_1))}
  =\|\iota_1(x_{0})-\iota_1(x_{pq})\|
  \le \sum_{r=1}^p
    \|\iota_1(x_{(r-1)q})-\iota_1(x_{rq})\|
  =p\sqrt{2(1-u_q(\theta_1))}.
\]
Dividing by $pq$ proves the first assertion.

Now fix $2\le i\le d$. Taking $p=i$ and $q=1$ gives
\[
  \frac{\sqrt{2(1-u_i(\theta_1))}}{i}
  \le \sqrt{2(1-u_1(\theta_1))}.
\]
Equality would imply that
$\iota_1(x_0),\ldots,\iota_1(x_i)$ lie on a line.
These points are distinct because $\iota_1$ is injective.
Since $i\ge2$, there are at least three such points. However, a line
intersects the unit sphere in at most two points.
Hence the inequality is strict.
\end{proof}

By~\eqref{eq:standard-distances},
the expansion ratio of $\iota_1$ on a pair at distance $i$ is
$\sqrt{2(1-u_i(\theta_1))}/i$.
Thus Lemma~\ref{lem:distance-range} shows that the maximum expansion
ratio is $\sqrt{2(1-u_1(\theta_1))}$, attained only on edges.
For $1\le i\le\lfloor d/2\rfloor$, put
$r=\lfloor d/i\rfloor$. Then $d/2<r i\le d$.
The first assertion of Lemma~\ref{lem:distance-range}, with $p=r$
and $q=i$, gives
\[
  \frac{\sqrt{2(1-u_{r i}(\theta_1))}}{r i}
  \le \frac{\sqrt{2(1-u_i(\theta_1))}}{i}.
\]
Thus the minimum over all distances is attained at a distance greater
than $d/2$, as in~\cite[Corollary~2.2]{CGIK}.
Consequently,
\[
  \min_{1\le i\le d}
    \frac{\sqrt{2(1-u_i(\theta_1))}}{i}
  =\min_{\lfloor d/2\rfloor+1\le i\le d}
    \frac{\sqrt{2(1-u_i(\theta_1))}}{i}.
\]
Since $\iota_1$ is injective, the definition of distortion and the
preceding observations give
\begin{align}
  (\dist \iota_1)^2
  &=\left(
      \max_{x\ne y}
      \frac{\|\iota_1(x)-\iota_1(y)\|}{\partial(x,y)}
    \right)^2
    \left(
      \max_{x\ne y}
      \frac{\partial(x,y)}{\|\iota_1(x)-\iota_1(y)\|}
    \right)^2 \notag\\
  &=2(1-u_1(\theta_1))
    \max_{\lfloor d/2\rfloor+1\le i\le d}
      \frac{i^2}{2(1-u_i(\theta_1))} \notag\\
  &=\max_{\lfloor d/2\rfloor+1\le i\le d}
      i^2\frac{1-u_1(\theta_1)}{1-u_i(\theta_1)}
  =\max_{\lfloor d/2\rfloor+1\le i\le d}
      \frac{i^2}{S_i(\theta_1)}.
  \label{eq:standard-distortion}
\end{align}

For a map $F:X\to\R^m$, define the mean squared distance
\[
  W_i(F)=\frac{1}{|X|k_i}
    \sum_{\substack{x,y\in X\\\partial(x,y)=i}}
      \|F(x)-F(y)\|^2\qquad(1\le i\le d).
\]
The sum is over ordered pairs.
For the standard representation $\iota_1$, we have
$W_i(\iota_1)=2(1-u_i(\theta_1))$.
The following Wirtinger inequality for the distance relations of
$\Gamma$ follows from Theorem~\ref{thm:central}.

\begin{proposition}
\label{lem:wirtinger}
For every nonconstant map $F:X\to\R^m$ and every $1\le i\le d$,
we have 
\[
  \frac{W_i(F)}{W_1(F)}\le
  S_i(\theta_1)=\frac{W_i(\iota_1)}{W_1(\iota_1)}.
\]
Moreover, suppose that $2\le i\le d$.
Let $G$ be the Gram matrix of $F$. Then equality holds if and only if $E_jG=0$ for all $2\le j\le d$.
\end{proposition}

\begin{proof}
Let $G=(\ip{F(x)}{F(y)})_{x,y\in X}$ be the Gram matrix of $F$.
Both $G$ and $E_j$ are positive semidefinite, so
$\operatorname{tr}(E_jG)\ge0$.
Using~\eqref{eq:idempotents} and the row sums of $A_i$, we obtain
\begin{align*}
  W_i(F)
  &=\frac{2}{|X|}\operatorname{tr}
       \left(\left(I-\frac{A_i}{k_i}\right)G\right)\\
  &=\frac{2}{|X|}\operatorname{tr}
       \left(\left(I-\sum_{j=0}^d u_i(\theta_j)E_j \right)G\right)\\
  &=\frac{2}{|X|}\sum_{j=1}^d
       (1-u_i(\theta_j))\operatorname{tr}(E_jG).
\end{align*}
Theorem~\ref{thm:central} gives
$1-u_i(\theta_j)\le S_i(\theta_1)(1-u_1(\theta_j))$
for $1\le j\le d$.
Therefore
\[
  W_i(F)\le\frac{2S_i(\theta_1)}{|X|}
    \sum_{j=1}^d(1-u_1(\theta_j))\operatorname{tr}(E_jG)
  =S_i(\theta_1)W_1(F).
\]
Since $\Gamma$ is connected and $F$ is nonconstant, $W_1(F)>0$.
Dividing by $W_1(F)$ proves the inequality.

Suppose now that $2\le i\le d$.
The equality condition in Theorem~\ref{thm:central} shows that equality
holds if and only if $\operatorname{tr}(E_jG)=0$ for every $2\le j\le d$.
Since $G$ is positive semidefinite and $E_j$ is an orthogonal projection,
this is equivalent to $E_jG=0$ for every $2\le j\le d$.
\end{proof}

We now combine the distortion of $\iota_1$ with the Wirtinger
inequality to prove Corollary~\ref{cor:distortion}.

\begin{proof}[Proof of Corollary~\ref{cor:distortion}]
By~\eqref{eq:standard-distortion}, the squared distortion of
$\iota_1$ is the right-hand side of~\eqref{eq:distortion}. This gives
the required upper bound for $c_2(\Gamma)^2$.

For the reverse bound, fix $1\le i\le d$.
Lemma~\ref{lem:monotone} gives $S_i(\theta_1)>0$.
Let $F$ be an injective map with distortion $D$. Rescale $F$ so that
\[
  \partial(x,y)\le\|F(x)-F(y)\|\le D\partial(x,y)
  \qquad(x,y\in X).
\]
Then $W_i(F)\ge i^2$ and $W_1(F)\le D^2$.
Proposition~\ref{lem:wirtinger} now gives
\[
  i^2\le W_i(F)\le S_i(\theta_1)W_1(F)
  \le S_i(\theta_1)D^2.
\]
Taking the infimum over $F$ and then the maximum over
$\lfloor d/2\rfloor+1\le i\le d$ gives the reverse bound.
Therefore~\eqref{eq:distortion} holds, and $\iota_1$ is optimal.
\end{proof}

\begin{remark}
The intersection array determines $k_i$, the eigenvalues, and the
cosine sequences. Corollary~\ref{cor:distortion} therefore expresses
$c_2(\Gamma)$ entirely in terms of the intersection array.
The proof uses no transitivity assumption on the automorphism
group of $\Gamma$.
\end{remark}

We conclude this section with an example of an optimal embedding
that is not similar to the standard representation in
Corollary~\ref{cor:distortion}.

We use the following construction of the Sylvester graph from
$\mathrm{PG}(2,9)$; see~\cite[Section~13.1]{BCN}.
Write $\mathbb F_9=\mathbb F_3(\iota)$, where $\iota^2=-1$.
A point of $\mathrm{PG}(2,9)$ is a one-dimensional subspace $[x]$
represented by a nonzero vector $x\in\mathbb F_9^3$.
Thus $x$ and $\lambda x$ represent the same point for every
$\lambda\in\mathbb F_9^\times$.
For $x=(x_1,x_2,x_3)^T$ and $y=(y_1,y_2,y_3)^T$, put
\[
  Q(x)=x_1^2+x_2^2+x_3^2,
  \qquad
  x\cdot y=x_1y_1+x_2y_2+x_3y_3.
\]
The equation $Q(x)=0$ defines a nondegenerate conic $\mathcal C$. Put
\[
  \mathcal I:=
  \{[x]\in\mathrm{PG}(2,9):
    Q(x)\in\mathbb F_9^\times\setminus(\mathbb F_9^\times)^2\}.
\]
Here $(\mathbb F_9^\times)^2$ is the set of nonzero squares in
$\mathbb F_9$.
This set is well defined, since
$Q(\lambda x)=\lambda^2Q(x)$.
A standard count gives $|\mathcal I|=36$.
Define a graph $\Gamma$ with vertex set $\mathcal I$.
Two distinct vertices $[x]$ and $[y]$ are adjacent when $x\cdot y=0$.
Equivalently, $[y]$ lies on the polar line
\[
  x^\perp=\{[z]\in\mathrm{PG}(2,9):x\cdot z=0\}
\]
of $[x]$ with respect to $\mathcal C$.
The adjacency condition is independent of the representatives.
It is symmetric because $x\cdot y=y\cdot x$.
Moreover, $x\cdot x=Q(x)\ne0$, so the graph has no loops.

Then this graph is distance-regular of diameter three, with
intersection array and spectrum
\[
  \{5,4,2;1,1,4\},\qquad
  5^{(1)},\quad 2^{(16)},\quad (-1)^{(10)},\quad (-3)^{(9)}.
\]
It is called the \emph{Sylvester graph}.

\begin{proposition}
\label{prop:sylvester}
The Sylvester graph has least Euclidean distortion $3/\sqrt2$.
It admits an optimal embedding whose affine span has dimension eight.
Its standard representation associated with $\theta_1$ has
affine span of dimension $16$.
\end{proposition}

\begin{proof}
The cosine sequence for $\theta_1=2$ is
$(1,2/5,-1/20,-1/5)$.
Thus $S_2(\theta_1)=7/4$ and $S_3(\theta_1)=2$.
Corollary~\ref{cor:distortion} gives
\[
  c_2(\Gamma)^2=\max\left\{\frac{16}{7},\frac92\right\}
  =\frac92.
\]
By~\eqref{eq:idempotents},
\begin{equation}
\label{eq:sylvester-E1}
  E_1=\frac1{45}(20I+8A_1-A_2-4A_3).
\end{equation}
The standard representation $\iota_1$ has Gram matrix
$(9/4)E_1$ and affine span of dimension $16$.

For vertices $[x],[y]\in\mathcal I$, choose vector representatives $x,y$
and define
\[
  \tau([x],[y])=\frac{(x\cdot y)^2}{Q(x)Q(y)}.
\]
This value is independent of the representatives.
A direct calculation shows that $\tau([x],[y])$ takes only the values
$1,0,\pm\iota,-1$.
These values occur precisely when $\partial([x],[y])=0,1,2,3$, respectively.
Define a real symmetric matrix $T$ by
\[
  T_{[x],[y]}=\begin{cases}
    1,&\tau([x],[y])=\iota,\\
    -1,&\tau([x],[y])=-\iota,\\
    0,&\text{otherwise}.
  \end{cases}
\]

We compute $T^2$ by counting projective points.
For $\varepsilon,\delta\in\{1,-1\}$, let
\[
  N_{\varepsilon\delta}([x],[y])=
  |\{[z]\in\mathcal I:\tau([x],[z])=\varepsilon\iota,
          \ \tau([z],[y])=\delta\iota\}|.
\]
The orthogonal group is transitive on ordered pairs with a fixed
value of $\tau$.
To see this, scale the representatives to have a fixed norm and
choose their signs. An isometry between their spans then extends to
the whole space.
We may therefore fix the vertex $[x]\in\mathcal I$ with representative
$x=(0,1,1+\iota)^T\in\mathbb F_9^3$.
Put $\nu=Q(x)=1-\iota$.
Each vertex has exactly two vector representatives $z$ with $Q(z)=\nu$.
Indeed, if $z_0$ represents a vertex, then
$Q(\lambda z_0)=\nu$ is equivalent to
$\lambda^2=\nu/Q(z_0)$. The quotient of the two nonsquares is a
nonzero square in $\mathbb F_9$, and hence this equation has two
solutions for $\lambda$.
For any representative $y$ of a vertex $[y]\in\mathcal I$,
$N_{\varepsilon\delta}([x],[y])$ is half the number of
solutions $z\in\mathbb F_9^3$ of
\[
  Q(z)=\nu,\qquad
  (x\cdot z)^2=\varepsilon\iota \nu^2= -\varepsilon,\qquad
  (y\cdot z)^2=\delta\iota Q(y)\nu.
\]
Counting the solutions gives the following table, where the signs in
the subscripts denote $1$ and $-1$.
\[
\begin{array}{c|c|rrr}
  \tau([x],[y])&y&N_{++}&N_{+-}=N_{-+}&N_{--}\\ \hline
  1&(0,1,1+\iota)^T&10&0&10\\
  0&(0,1,1-\iota)^T&4&0&4\\
  \iota&(1,0,1+\iota)^T&2&3&3\\
  -\iota&(1,1,-1-\iota)^T&3&3&2\\
  -1&(0,1,-1+\iota)^T&2&4&2
\end{array}
\]

As an example, we explain the count for the row $\tau([x],[y])=0$.
Take $y=(0,1,1-\iota)^T$ and write $z=(z_1,z_2,z_3)^T$.
Put
\[
  r=x\cdot z=z_2+(1+\iota)z_3,\qquad s=y\cdot z=z_2+(1-\iota)z_3.
\]
Since $Q(y)=1+\iota$, the two inner product
conditions become
\[
  r^2=-\varepsilon,\qquad s^2=-\delta\iota.
\]
Each equation has two solutions in $\mathbb F_9$.
Each pair $(r,s)$ determines $z_2,z_3$ uniquely.
Since $Q(z)=z_1^2+r^2/\nu+s^2/(1+\iota)$, the condition $Q(z)=\nu$
becomes
\[
  z_1^2=\nu\bigl(1-(\varepsilon+\delta)\iota\bigr).
\]
For $(\varepsilon,\delta)=(1,1)$ or $(-1,-1)$, the right-hand side
is $-1$ or $\iota$, respectively.
Both are nonzero squares, so each case gives $4\cdot2=8$ vectors.
If $\varepsilon=-\delta$, the right-hand side is the nonsquare $\nu$,
so there are no solutions.
Dividing these counts by two gives $N_{++}=N_{--}=4$ and
$N_{+-}=N_{-+}=0$.

Since $(T^2)_{[x],[y]}=N_{++}+N_{--}-N_{+-}-N_{-+}$,
the table and~\eqref{eq:sylvester-E1} yield
\begin{equation}
\label{eq:sylvester-square}
  T^2=20I+8A_1-A_2-4A_3=45E_1.
\end{equation}

The matrix $T$ is real symmetric and has zero diagonal.
Equation~\eqref{eq:sylvester-square} shows that it vanishes on
the orthogonal complement of $E_1\R^X$.
On $E_1\R^X$, its eigenvalues are $3\sqrt5$ and $-3\sqrt5$.
Each has multiplicity eight, since $\operatorname{tr}T=0$.
It follows that
\begin{equation}
\label{eq:sylvester-Gram}
  G=\frac{15}{8}E_1+\frac{\sqrt5}{8}T
\end{equation}
is positive semidefinite of rank eight.
Moreover, $G\1=0$.
Let $F:X\to\R^8$ have Gram matrix $G$.
Its centroid is zero, and its affine span is $\R^8$.
Since $T=E_1TE_1$, we have $E_2G=E_3G=0$.
Thus $F$ attains equality in Proposition~\ref{lem:wirtinger} for $i=2,3$.
For distinct vertices $[x],[y]$, its squared distances are
\[
  \|F([x])-F([y])\|^2=
  \begin{cases}
    1,&\partial([x],[y])=1,\\
    \dfrac74-\dfrac{\sqrt5}{4}T_{[x],[y]},&\partial([x],[y])=2,\\[2pt]
    2,&\partial([x],[y])=3.
  \end{cases}
\]
For pairs at distance two, we have
\[
  \frac29<\frac{7-\sqrt5}{16}
  \le\frac{\|F([x])-F([y])\|^2}{4}
  \le\frac{7+\sqrt5}{16}<1.
\]
Thus $F$ is injective.
Its maximum expansion ratio is one, attained on every edge.
Its maximum contraction ratio is $3/\sqrt2$, attained at distance three.
Hence $F$ is optimal.
The dimensions of the affine spans show that $F$ and $\iota_1$ are
not similar.
\end{proof}

\begin{remark}
Diameter three is the smallest possible for this nonuniqueness.
For diameter two, let $F$ be optimal and scale it so that its
maximum expansion ratio is one.
Corollary~\ref{cor:distortion} and the Wirtinger inequality give
\[
  S_2(\theta_1)=\frac4{c_2(\Gamma)^2}
  \le W_2(F)\le S_2(\theta_1)W_1(F)\le S_2(\theta_1).
\]
Equality forces every edge to have squared length one and every
pair at distance two to have squared length $S_2(\theta_1)$.
These are all the pairs, so $F$ is similar to $\iota_1$.
For diameter one, every optimal embedding is a regular simplex
up to similarity.
\end{remark}

\section*{Acknowledgments}
The author is supported by JSPS KAKENHI Grant Number JP24K06830.

\bibliographystyle{plain}
\bibliography{ref}

\end{document}